\documentclass[a4paper,12pt,fleqn]{article}

\usepackage{theorem}
\usepackage{amssymb}
\usepackage{latexsym}

\theoremstyle{break}
\newtheorem{thm}{Theorem}[section]
\newtheorem{prop}{Proposition}[section]
\newtheorem{lem}{Lemma}[section]
\newtheorem{de}{Definition}[section]

\newtheorem{nota}{Notation}[section]

\title{ NBB bases of some pattern avoiding lattices
 }
\author{Masaya Tomie}
\date{Institute of Mathematics, University of Tsukuba, 
Tsukuba, Ibaraki 305-8571, Japan (e-mail: tomie@math.tsukuba.ac.jp)}

\begin{document}

\maketitle

\begin{center}
{\bf Abstract}
\end{center}
In this paper we will determine the NBB bases with respect to 
standard ordering of coatoms (resp.atoms) of 
123-132-213-avoiding (resp.321-avoiding) lattices. 
Using these expression we will calculate the M\"obius numbers of  
123-132-213-avoiding lattices and 321-avoiding lattices. 
These values become some modification of fibonacci polynomials.

\section{Introduction}

In this paper we give expressions of NBB bases 
of 123-132-213-avoiding lattices and 321-avoiding lattices.
 Using these expressions we will calculate the  M\"obius numbers 
 of these lattices. Surprisingly these values become some 
 modification of fibonacci polynomials.
 We introduce modified fibonacci polynomials
$\{ F_{n}(q) \}_{n \in \mathbb{N}}$ later.
Let  $A_{n}$ (resp.$B_{n}$) be  the partially ordered set of the 
 123-132-213-avoiding (resp.321-avoiding) permutations with the weak 
  order on the permutation group 
  $S_{n}$ with a unique minimal element (resp.maximal element) 
  for $n \in \mathbb{N}$. 
We also determine the NBB bases for  $A_{n}$ and $B_{n}$ with 
respect to a natural total ordering of atoms or coatoms of them.
Using the modified fibonacci polynomials and the expression of the 
NBB bases 
 we will consider the M\"obius numbers of $A_{n}$ and $B_{n}$ for each 
 $n \in \mathbb{N}$.
 
  Let $P$ be a poset and  ${\rm Int}(P)$ the set of intervals of $P$.
   We call the function 
  $\mu : \ {\rm Int}(P) \rightarrow \mathbb{Z}$ the M\"obius function of $P$ 
  if $\mu$ satisfies the following identity.
  
  \begin{eqnarray}
  \sum_{x \le y \le z}\mu([x,y]) = \delta_{x,z}
  \end{eqnarray}
  
  If $P$ has  a  maximum element $\widehat{1}$ and a 
  minimum element $\widehat{0}$. Then we put 
  $\mu(P) := \mu([\widehat{0},\widehat{1}])$. We call $\mu(P)$ the M\"obius 
  number of $P$.
 Our main result is as following.

 \begin{thm}
 
 For $n \in \mathbb{N}_{\ge 3}$ we have 
 
 \begin{eqnarray}
 \mu(A_{n}) = \mu(B_{n}) = (-1) F_{n-2}(-1).
 \end{eqnarray}

 \end{thm}

\section{Preliminaries}

\subsection{Bounded below sets}

This subsection we introduce a technique to calculate 
 M\"obius numbers of lattices which is given in Blass and Sagan's paper 
\cite{blasssagan}.

Throughout this subsection $(L, \le )$  will 
denote a finite lattices. We will denote it $L$ for short. We will use 
$\wedge$ for the meet (greatest lower bound) and $\vee$ for the join 
(least upper bound) in $L$.
Since $L$ is finite it also has the unique minimal element $\widehat{0}$ and 
the unique maximal element $\widehat{1}$. 
We let $\mu(L) := \mu([\widehat{0},\widehat{1}])$.
Our goal in this subsection is to give a combinatorial description of 
$\mu(L)$.
Let $A(L)$ (resp.$B(L)$)
be the set of coatoms (resp.atoms) of $L$.
Give $A(L)$ (resp.$B(L)$) 
 an arbitrary total order, which we denote $\unlhd_{A}$ 
 (resp.$\unlhd_{B}$) to distinguish 
it from $\le$ in $L$.
 A nonempty set $D \subseteq A(L)$ (resp.$D^{\prime} \subseteq B(L)$) is 
 $bounded \ below$ ($BB$ for short) 
 if for 
 every $d \in D$ (resp.$d^{\prime} \in D^{\prime}$)
  there is an $a \in A(L)$ (resp. $a^{\prime} \in B(L)$)  
 such that 
 $a \lhd_{A} d$ and $a > \wedge D$ (resp.$a^{\prime}
  \lhd_{B} d$ and $a^{\prime} < \vee D^{\prime}$ ).
  So $a$ (resp.$a^{\prime}$) is simultaneously a strict lower bound 
  for $d$ (resp.$d^{\prime}$)
   in the total order $\unlhd_{A}$ (resp.$\unlhd_{B}$)
   and for $\wedge D$ (resp.$\vee D^{\prime}$)
  in $\le$.
  We will say that $B \subseteq A(L)$ (resp.$B^{\prime} \subset B(L)$) is 
  $NBB$ if $B$ (resp.$B^{\prime}$) does not contain any $D$ (resp.$D^{\prime}$)
  which is bounded below. In this case we will call $B$ (resp.$B^{\prime}$
  ) an $NBB \  base$ for 
  $x = \wedge B$ (resp.$x^{\prime} = \vee B^{\prime}$). In \cite{blasssagan} 
  Blass and Sagan proved the following statement.

  \begin{thm}[\cite{blasssagan}]\label{blass-sagan-main}
  
  Let $L$ be any finite lattice. 
  Let $A(L)$ (resp.$B(L)$ be the set of coatoms (resp.atoms) of $L$
  and $\unlhd_{A}$ (resp.$\unlhd_{B}$) any total order on 
  $A(L)$ (resp.$B(L)$). Then we have 
  
  \begin{eqnarray}
  \mu(L) &=& \sum_{B \in A(L), \ \wedge B = \widehat{0}} (-1)^{|B|}\\
  \mu(L) &=& \sum_{C \in B(L), \ \vee C = \widehat{1}} (-1)^{|C|}
  \end{eqnarray}
  
  where the sum is over all $NBB$ bases  of $\widehat{0}$ (resp.$\widehat{1}$)
   and $| \cdot |$ denotes cardinality.

  \end{thm}

\subsection{ Modified Fibonacci polynomials }

In this subsection  we introduce modified Fibonscci polynomials.

\begin{de}\label{def of q-fibonacci}

We define the sequences $\{ F_{n}(q) \}_{n \in \mathbb{N}}$ by the following 
relations:
\begin{eqnarray}
F_{1}(q) :&=& 1, \ F_{2}(q):= 1, \\
F_{k+2}(q) &=& F_{k+1}(q) + q F_{k}(q), \ {\rm for}\  k \ge 1.
\end{eqnarray}
\end{de}
We call the sequences $\{ F_{n} \}_{n \in \mathbb{N}}$  
modified Fibonacci polynomials.

\begin{nota}

The Fibonacci polynomials $F_{n}^{\prime}(q)$ are defined by the following 
relations:
\begin{eqnarray}
F^{\prime}_{1}(q) :&=& 1, \ F^{\prime}_{2}(q):= q, \\
F^{\prime}_{k+2}(q) &=& q  F^{\prime}_{k+1}(q) +  
F^{\prime}_{k}(q), \ {\rm for}\  k \ge 1.
\end{eqnarray}
But in this paper we don't use the Fibonacci polynomials.

\end{nota}

Let $X$ be a subset of $[n] := \{1,2, \ldots, n \}$. We call $X$ a sparse set 
if and only if $1 \in X$ and if $i \in X$ then  $i+1 \not\in X$ for 
$1 \le i \le n-1$.
 For $n=4$ the corresponding sparse sets are $\{1 \}, \ \{1, 3 \} $ and 
 $\{1,4 \}$.
 
Then we have the following proposition. 
A simple calculation yields the statement of the proposition so we omit 
the proof.

\begin{prop}
 
Let $H_{n}(q) := \sum_{X: {\rm sparse \ set \ of \ } [n] } q^{\sharp X -1}$.
Then we have 
\begin{eqnarray}
H_{n}(q) = F_{n}(q)
\end{eqnarray}
for $n \in \mathbb{N}$.

\end{prop}

\subsection{The weak order on the symmetric group}

In this subsection we will introduce the weak 
order and its lattice structure \cite{aguiar-sottile}
\cite{guilbaud-rosenstiehl}.
For $n \in \mathbb{N}$ let $\sigma$ be an element of the permutation group 
 $S_{n}$. We put 
${\rm Inv}(\sigma) := 
\{ (i,j) \ | \ 1 \le i < j \le n, \ \sigma(i) > \sigma(j) \}$.
 We write $\sigma \le \tau$  if ${\rm Inv}(\sigma) \subset {\rm Inv}(\tau)$.
 This determines the $weak  \ order $ on $S_{n}$.
 This weak order is a lattice. The identity permutation $1_{n}$ is the 
 minimum element and $n (n-1) \cdots 21$ is the maximum element. 
 A set $J$ is the inversion set of a permutation in $S_{n}$ if and only if 
 both $J$ and its complement ${\rm Inv} (n (n-1) \cdots 21) - {\rm Inv}(J)$ 
 are transitively closed (i.e. $(i,j) \in J$ and $(j,k) \in J$ imply 
 $(i,k) \in J$, and the same for its complement). 
 The join (least upper bound) of 
 two permutations $\sigma$ and $\tau \in S_{n}$ is the 
 permutation 
 whose inversion set is the transitive closure of the union of the inversion 
 sets of $\sigma$ and $\tau$ 
 \begin{eqnarray}
 \{ (i,j) \ | \ \exists {\rm chain} \ i= k_{0}< \cdots <k_{s} = j 
 \ s.t. \ \forall r, \ (k_{r-1},k_{r}) \in {\rm Inv} (\sigma) \cup 
 {\rm Inv} (\tau) \}.
 \end{eqnarray}
We denote it $\sigma \vee \tau$.

Similarly, the meet (greatest lower bound) of $\sigma$ and $\tau$ is the 
permutation whose inversion set if 
\begin{eqnarray}
\{ (i,j) \ | \ \forall  {\rm chains} \ i= k_{0} < \cdots k_{s} = j, \exists 
r \ s.t. \ (k_{r-1},k_{r} ) \in {\rm Inv}(\sigma) \cap {\rm Inv}(\tau) \}.
\end{eqnarray}

\section{The case of 123-132-213 avoiding lattices}

For each $n \in \mathbb{N}$ we define $A^{\prime}_{n}$ to be the 
partially ordered set of 123-132-213-avoiding permutations associated with 
the weak order on $S_{n}$. We put $A_{n} := A_{n}^{\prime} \cup \{
\widehat{0} \}$ where $\widehat{0}$ is a unique minimum element.
For example we have 
  $A_{1} = \{ 1, \widehat{0} \}$,  $A_{2} = \{ 12, 21, \widehat{0} \}$ and 
  $A_{3} = \{ 231, 312, 321, \widehat{0} \}$.
Let $X$ be a subset of permutations of $S_{n}$. 
Put $(n+1)X := \{ (n+1) a_{1} \cdots a_{n} \ | \ a_{1} \cdots a_{n} \in 
X \ \}$. 
The following theorem is known.

\begin{thm}[\cite{barcucci-bernini-poneti}]\label{type-A-shape}

Let $\tilde{A_{n}}$ be the set of 123-132-213-avoiding  permutations in 
$S_{n}$ for $n \in \mathbb{N}$. 
(We will consider $\tilde{A_{n}}$ as a set.)
 Then we have 
\begin{eqnarray}
\tilde{A}_{n+2} = (n+1)(n+2)\tilde{A}_{n} \uplus (n+2)\tilde{A}_{n+1}
\end{eqnarray}

\end{thm}

From Theorem {\rmfamily \ref{type-A-shape}} we have the following lemma.

\begin{lem}
The poset $A_{n}$ is an order filter of $S_{n}$ for $n \in \mathbb{N}$.
\end{lem}

{\bf PROOF}

We will prove by induction on $n$. We assume that this lemma holds for 
$\le n$.
Let $\sigma \in A_{n+1}$. 
Then we have $\sigma = n (n+1) a_{1} \dots a_{n-1}$ with 
$a_{1} \dots a_{n-1} \in A_{n-1}$
 or 
 $\sigma = (n+1) b_{1} \cdots b_{n}$ with $b_{1} \cdots b_{n} \in A_{n}$
  by Theorem {\rmfamily \ref{type-A-shape}}.

The case of  $\sigma = n (n+1) a_{1} \dots a_{n-1}$ with 
$a_{1} \dots a_{n-1} \in A_{n-1}$. For $\tau$ with $\tau \ge \sigma$ in the 
weak order on $S_{n+1}$, 
we have 
either $\tau =  n (n+1) a^{\prime}_{1} \cdots a^{\prime}_{n-1}$  or 
$\tau =  (n+1) n  a^{\prime \prime}_{1} \cdots a^{\prime \prime}_{n-1}$ 
with $a^{\prime}_{1} \cdots a^{\prime}_{n-1}, 
a^{\prime \prime}_{1} \cdots a^{\prime \prime}_{n-1} 
 \ge a_{1} \dots a_{n-1}$.
By assumption we have $ a^{\prime}_{1} \cdots a^{\prime}_{n-1} , \ 
a^{\prime \prime}_{1} \cdots a^{\prime \prime}_{n-1} \in A_{n-1}$. 
Hence we have  $\tau \in A_{n+1}$.

The case of 
$\sigma = (n+1) b_{1} \cdots b_{n}$ with $b_{1} \cdots b_{n} \in A_{n}$.
For $\tau$ with   $\tau \ge \sigma$ in the 
weak order on $S_{n+1}$  we have 
$\tau = (n+1) b^{\prime}_{1} \cdots b^{\prime}_{n}$ with 
$b^{\prime}_{1} \cdots b^{\prime}_{n} \ge b_{1} \cdots b_{n} $.
By assumption we have $b^{\prime}_{1} \cdots b^{\prime}_{n} \in A_{n}$.
Hence we have $\tau \in A_{n+1}$.

This completes the proof of our lemma.
\ \ \ \ \ \  \  \ $\Box$

Next we define 
\begin{eqnarray}
C_{n} := \{ c_{1} c_{2} \cdots c_{k} \ | \ 
c_{i} = 1 \ {\rm or } \ 2 \ {\rm with } \ c_{1} + c_{2} + \cdots c_{k} = n \}
 \cup \{ \widehat{0} \}
\end{eqnarray}
with covering relations as following;
\begin{eqnarray}
 c_{1}  \cdots c_{i-1} 2 c_{i+1} \cdots c_{k} \prec 
 c_{1}  \cdots c_{i-1} 1 1  c_{i+1} \cdots c_{k}, \ \ 
\end{eqnarray}
where $\widehat{0}$ the minimum elements for  $n \in \mathbb{N}$.
Then $C_{n}$ has a poset structure for $n \in \mathbb{N}$.

\begin{prop}
For $n \in \mathbb{N}$
we have $A_{n} \simeq C_{n}$ as a poset.

\end{prop}

{\bf Proof}

For each $n \in \mathbb{N}$ we will define the map 
$\phi_{n}: \ C_{n} \rightarrow A_{n}$  by induction. 

We put $\phi_{1}(1) := 1, \ \phi_{1}(\widehat{0}) := \widehat{0}, \ 
\phi_{2}(11) := 21, \ \phi_{2}(2):= 12$ and  
$\phi_{2}(\widehat{0}) := \widehat{0}$. 

For $n \in \mathbb{N}, \ c_{1} \cdots c_{k} \in C_{n}$ we define 
$\phi_{n}$ as follows;
\begin{eqnarray}
\phi_{n}(c_{1} \cdots c_{k}) := n \phi_{n-1}(c_{2} \cdots c_{k}) \ \ \ \ \ 
{\rm if \ } c_{1} = 1  \\
\phi_{n}(c_{1} \cdots c_{k}) := (n-1) n 
 \phi_{n-2}(c_{2} \cdots c_{k}) \ \ \ \ \ 
{\rm if \ } c_{1} = 2 \\
\phi_{n}(\widehat{0}) := \widehat{0}.
\end{eqnarray}
For example we have $\phi_{3}(\widehat{0}) = \widehat{0}, \ 
\phi_{3}(111) = 321, \ \phi_{3}(21) = 231$ and $ \phi_{3}(12) = 312 $.

Next we will define the map 
$\psi_{n}: \ A_{n} \rightarrow C_{n}$ for each 
$n \in \mathbb{N}$ by induction. 

We put $\psi_{1}(1) := 1, \ \psi_{1}(\widehat{0}) := \widehat{0}, \ 
\psi_{2}(21) := 11, \ \psi_{2}(12):= 2$ 
and $ \psi_{2}(\widehat{0}) := \widehat{0}$. 

For $n \in \mathbb{N}$ we define $\psi_{n}$ as follows;
\begin{eqnarray}
\psi_{n}((n-1) n a_{1} \cdots a_{n-2}) := 2 
\psi_{n-2}(a_{1} \cdots a_{n-2})   \\
\psi_{n}(n a_{1} \cdots a_{n-1}) := 1 
 \psi_{n-1}(a_{1} \cdots a_{n-1}) \\
\psi_{n}(\widehat{0}) := \widehat{0}.
\end{eqnarray}
For example we have $\psi_{3}(\widehat{0}) = \widehat{0}, \ 
\psi_{3}(321) = 111, \ \psi_{3}(231) = 21$ and $  \psi_{3}(312) = 12$.
Note that if $a_{1} \cdots a_{n} \in A_{n}$ we have either $a_{1} = n-1$ and 
$a_{2} = n$  or $a_{1} = n$.
It is easy to see that $\psi_{n} \circ \phi_{n} = id_{C_{n}}$ and 
$\phi_{n} \circ \psi_{n} = id_{A_{n}}$.

Next we have to show that the map $\phi_{n}$ and $\psi_{n}$ 
are both order preserving.

The case of $\phi_{n}$. We will show that $\phi_{n}$ 
preserves the covering relation of $C_{n}$.
For $c_{1} \cdots c_{i-1} 2 c_{i+1} \cdots c_{k}, \ 
c_{1} \cdots c_{i-1} 1 1 c_{i+1} \cdots c_{k} \in C_{n}$ we have 
$c_{1} \cdots c_{i-1} 2 c_{i+1} \cdots c_{k} 
\prec c_{1} \cdots c_{i-1} 1 1 c_{i+1} \cdots c_{k}$. 
We have $\phi_{n}(c_{1} \cdots c_{i-1} 1 1 c_{i+1} \cdots c_{k}) = 
\phi_{n-2-x}^{x}(c_{1} \cdots c_{i-1}) 
(x+2) (x+1) \phi_{n}(c_{i+1} \cdots c_{k})$
 and 
 $\phi_{n}(c_{1} \cdots c_{i-1} 2 c_{i+1} \cdots c_{k}) = 
\phi_{n-2-x}^{x}(c_{1} \cdots c_{i-1}) 
(x+1) (x+2) \phi_{n}(c_{i+1} \cdots c_{k})$ where 
$x = c_{i+1} + \cdots + c_{k}$ and $\phi_{n-2-x}^{x}(c_{1} \cdots c_{i-1})
= (c^{\prime}_{1} + x) \cdots (c^{\prime}_{i-1} + x)$ for 
$\phi_{n-2-x}(c_{1} \cdots c_{i-1})
= c^{\prime}_{1} \cdots c^{\prime}_{i-1}$.
It is easy to see that 
$\phi_{n-2-x}^{x}(c_{1} \cdots c_{i-1}) 
(x+1) (x+2) \phi_{n}(c_{i+1} \cdots c_{k}) \prec 
\phi_{n-2-x}^{x}(c_{1} \cdots c_{i-1}) 
(x+2) (x+1) \phi_{n}(c_{i+1} \cdots c_{k})$ in $A_{n}$.

The case of $\psi_{n}$. We will show that $\psi_{n}$ 
preserves the covering relation of $A_{n}$.
Let $a_{1} \cdots a_{n} \in A_{n}$ with $a_{i} = m$ and $a_{j} = m+1$ for 
$i < j$.  Because $a_{1} \cdots a_{n} $ avoids 213 pattern and 132 pattern, 
we have $j = i+1$.
We have $a_{1} \cdots a_{i-1} m (m+1) a_{i+2} \cdots a_{n} \prec 
a_{1} \cdots a_{i-1}  (m+1) m  a_{i+2} \cdots a_{n}$. Then we have  
$\phi_{n}(a_{1} \cdots  a_{i-1} m (m+1) a_{i+2}   \cdots a_{n} ) = 
\phi_{i-1}( st (a_{1} \cdots  a_{i-1}) ) 2 
\phi_{n-1-x} (st(a_{i+2} \cdots a_{n}))$ and 
$\phi_{n}(a_{1} \cdots  a_{i-1}  (m+1) m a_{i+2}   \cdots a_{n} ) = 
\phi_{i-1}( st (a_{1} \cdots  a_{i-1}) ) 1 1 
\phi_{n-1-i} (st(a_{i+2} \cdots a_{n}))$ where 
$ st(a_{1} \cdots  a_{i-1}) \in S_{i-1}$ is the unique permutation 
$\sigma \in S_{i-1}$ such that $\sigma_{s} < \sigma_{t} 
\Leftrightarrow a_{s} < a_{t}$.
It is easy to see that $\phi_{i-1}( st (a_{1} \cdots  a_{i-1}) ) 2 
\phi_{n-1-x} (st(a_{i+2} \cdots a_{n})) \prec 
\phi_{i-1}( st (a_{1} \cdots  a_{i-1}) ) 1 1 
\phi_{n-1-x} (st(a_{i+2} \cdots a_{n}))$.

This completes the proof of our proposition. \ \ \ \ \ \ \ $\Box$

Next we calculate the M\"obius numbers of $C_{n}$ for each $n \in \mathbb{N}$.

Let $A(C_{n}) := \{ 21 \cdots 1, \ 121 \cdots 1, \ldots 1 \cdots 21, 
1 \cdots 12 \} $ be the set of coatoms of $C_{n}$ 
where $\sharp A(C_{n}) = n-1$. We give $A(C_{n})$ a tatal order 
$\lhd$  as follows,
\begin{eqnarray}
21 \cdots 1 \lhd \ 121 \cdots 1 \lhd  \ldots  \lhd 1 \cdots 21 \lhd  
1 \cdots 12.
\end{eqnarray}
We put $\theta_{i} := 1 \cdots 1 
\underbrace{2}_{i-{\rm th}} 1 \cdots 1$. Then the following lemma is easy 
to prove so we will omit the proof.

\begin{lem}

For $1 \le i \le n-2$ we have 
\begin{eqnarray}
\theta_{i} \wedge \theta_{i+1} = \widehat{0}.
\end{eqnarray}

\end{lem}

\begin{lem}

We let $\{  i_{1}, \ldots i_{k} \} \subset [n-1]$ with $i_{1} < \ldots 
i_{k}$. If $i_{p+1} > i_{p} + 1$ for $1 \le p \le k-1$ then we have 
\begin{eqnarray}
\theta_{i_{1}} \wedge \theta_{i_{2}} \wedge \cdots \wedge \theta_{i_{k}} 
= 1 \cdots  \underbrace{2}_{j_{1}}  \cdots  
\underbrace{2}_{j_{2}}  \cdots \underbrace{2}_{j_{k}} \cdots 1,
\end{eqnarray}
where $j_{1} = i_{1}, \ j_{2} = i_{2}-1, \ j_{3} = i_{3} -2, \ 
\ldots, j_{k} = i_{k} - k +1$.
\end{lem}

{\bf Proof}

For $j \ge i+2$ we have 

$a_{1} \cdots a_{i-1}  \underbrace{2}_{i-{\rm th}} 1 \cdots 1 \wedge 
a_{1} \cdots a_{i-1} 1  \cdots \underbrace{2}_{j-{\rm th}} \cdots 1 = 
a_{1} \cdots a_{i-1} \cdots
 \underbrace{2}_{i-{\rm th}} \cdots \underbrace{2}_{(j-1){\rm th}} 
\cdots 1$ 

for $a_{1}, \ldots a_{i-1} \in \{ 1, \ 2 \}$.
From this fact and using induction we obtain the desired result.
 \ \ \ \ \ \ \ \ \ $\Box$

\begin{lem}\label{nbb-condition}

Put $X := \{ \theta_{i_{1}}, \ldots \theta_{i_{k}} \} \subset 
A(C_{n})$ with 
$i_{p+1} > i_{p} + 1$ for $1 \le p \le k-1$. Then the set $X$ is not BB 
 with respect to our total ordering $\lhd$.
 Moreover $X$ is NBB.

\end{lem}

{\bf Proof}

We have 
$\wedge X = 1 \cdots \underbrace{2}_{j_{1}} 1 \cdots \underbrace{2}_{j_{2}}
\cdots 1 \underbrace{2}_{j_{k}} \cdots 1$
where $j_{1} = i_{1}, j_{2} = i_{2} -1, \ldots ,j_{k} = i_{k} - k +1$. 
Then we have $\{ y \in A(C_{n}) \ | \ y \le \wedge X \} = 
\{ \theta_{i_{1}}, \ldots \theta_{i_{k}}  \}$. This yields that 
$X$ is not BB. For $Y \subset X$ the same argument yelds that $Y$ is not 
$BB$.
Hence we obtain the derived result. \ \ \ \ \ \ \ $\Box$

\begin{thm}\label{123-132-213-main-prop}
Let 
$X := \{  \theta_{i_{1}}, \ldots \theta_{i_{k}} \} $ be a subset of 
$A(C_{n})$ where $i_{1} < \ldots  < i_{k}$.
Then $X$ is an NBB base of 
$\widehat{0}$ 

$\Longleftrightarrow$ $X$ satisfies 
\begin{enumerate}
\item $i_{1} = 1$, \\
\item $\{ i_{2}-1, i_{3}-1, \ldots i_{k}-1 \}$ is a sparse set in $[n-2]$.
\end{enumerate}

\end{thm}

{\bf Proof}

($\Longrightarrow$)

If $X := \{  \theta_{i_{1}}, \ldots \theta_{i_{k}} \} $ is an NBB base of 
$\widehat{0}$. For each $1 \le i \le n-1$ we have $a_{i} > \wedge X = 
\widehat{0}$.  So we have $\theta_{1} \in X $ because  $X$ is  BB.
On the other hand bacause $\theta_{i_{1}} \wedge \cdots \wedge \theta_{i_{k}}
= \widehat{0}$ there exists $1 \le j \le k-1$ such that 
$i_{j+1} = i_{j}+1$. Then we have $\theta_{i_{j}} \wedge \theta_{i_{j+1}} 
= \widehat{0}$ and $\{ \theta_{i_{j}} , \theta_{i_{j+1}} \} \subset X$ 
is not BB. Hence we have $j=1, \ i_{1} = 1$ and $i_{2} = 2$. 
If there exists $2 \le j^{\prime} \le k-1$ such that 
$i_{j^{\prime}+1} = i_{j^{\prime}} +1$. We put $Z := \{ 
\theta_{i_{j^{\prime}}} , \ \theta_{i_{j^{\prime}+1}} \}$. Then 
$\theta_{i_{j^{\prime}}} \wedge \ \theta_{i_{j^{\prime}+1}} = \widehat{0}$ 
and $Z$ is not BB. This yields that $\theta_{1} \in Z$. This contaradicts 
the assumption $j^{\prime} \ge 2$. This completes the proof of 
"$\Longrightarrow$" part.

($\Longleftarrow$)

Let $X$ be a subset of $[n-1]$ satisfying the above conditions. 
For any subset $Y$ of $X$ we will show that  $Y$ is not BB.
We put $Y:= \{ \theta_{j_{1}}, \ldots , \theta_{j_{l}} \} \subset X$.

When $j_{1} = 1$ and $j_{2} = 2$ we have $\wedge Y = \widehat{0}$ and 
$\theta_{1} \in Y$. Hence $Y$ is not BB.
When $j_{1} = 1$ and $j_{2} \ge 3$ it is easy to see that 
$Y$ is not BB. 
When $j_{1} \ge 2$ 
it is also trivial from Lemma {\rmfamily \ref{nbb-condition}}. 
This completes the proof of our stetement.  \ \ \ \ \ \ \ $\Box$

From Theorem {\rmfamily \ref{blass-sagan-main}} and Theorem 
{\rmfamily \ref{123-132-213-main-prop}} we obtain the following result.

\begin{thm}

We have 

\begin{eqnarray}
\mu(A_{n}) = \mu(C_{n}) = \sum_{X \subset [n-2] \ X: \ {\rm sparse \ set}}
(-1)^{|X|+1} = (-1) F_{n-2}(-1).
\end{eqnarray}

\end{thm}

\section{The case of 321-avoiding lattices}

For each $n \in \mathbb{N}$ we define $B^{\prime}_{n}$ to be the 
partially ordered set of 321 avoiding permutations associated with 
the weak order on $S_{n}$. We put $B_{n} := B_{n}^{\prime} \cup \{
\widehat{1} \}$ where $\widehat{1}$ is a unique maximum element.
For example we have $B_{1} = \{ 1, \widehat{1} \} , \ 
 B_{2} = \{ 12, 21, \widehat{1} \}$ and $
 B_{3} = \{ 123,  213, 132, 312, 231, \widehat{1} \}$.
Lemma {\rmfamily \ref{321-order-ideal}} and Lemma 
{\rmfamily \ref{321-permutation}} 
are trivial from the definition of $B_{n}$.

\begin{lem}\label{321-order-ideal}

For each $n \in \mathbb{N}$ our poset $B^{\prime}_{n}$ is an order 
ideal of $S_{n}$. Therefore $B_{n}$ is a lattice.

\end{lem}

\begin{lem}\label{321-permutation}

For each $n \in \mathbb{N}$ we have 
\begin{eqnarray}
\sigma \in B^{\prime}_{n} 
\Longleftrightarrow 
{\rm If } \ (i,j) \in {\rm Inv}(\sigma) \ {\rm then} \ 
{\rm for \ any } \   k \ {\rm we \ have } \ 
(j,k) \not\in {\rm Inv}(\sigma).
\end{eqnarray}

\end{lem}

We put $\sigma_{i} := (i,i+1)$. 
Let 
$A(B_{n}) $ be
 the set of atoms of $B_{n}$. Note that 
 $A(B_{n})= \{ \sigma_{1}, \sigma_{2}, \ldots \sigma_{n-1} \}$.
We define $A(B_{n})$ a total order $\lhd$ as following;
\begin{eqnarray}
\sigma_{1} \lhd \sigma_{2} \lhd \ldots \lhd \sigma_{n-1}.
\end{eqnarray}
\begin{lem}\label{321-avoiding-neighbor}

We have 
$\sigma_{i} \vee \sigma_{i+1} = \widehat{1}$ for $1 \le i \le n-2$.

\end{lem}

{\bf Proof}

We fix $1 \le i \le n-2$.
We put $\tau := \sigma_{i} \vee \sigma_{i+1}$ in $S_{n}$ with the weak 
order. Then we have $\tau(i) > \tau(i+1) > \tau(i+2)$. 
 So we have $\tau \not\in 
B^{\prime}_{n}$. 
From Lemma {\rmfamily \ref{321-order-ideal}}
 we obtain the derived result. \ \ \ \ \ \ 
$\Box$

\begin{lem}\label{321-avoiding-sparse}
We assume   $1 \le j_{1} < \ldots < j_{l} \le n-1$ and 
$2 \le j_{p+1} - j_{p}$ for $1 \le p \le l-1$.
We put $Y := \{ \sigma_{j_{1}}, \ldots , \sigma_{j_{l}} \} \subset A(B_{n})$.
Then we have $\vee Y = 
\sigma_{j_{1}} \sigma_{j_{2}} \ldots  \sigma_{j_{l}}$.

{\bf Proof}

For $1 \le p \le q \le l$ we have $\sigma_{j_{p}} \sigma_{j_{q}} 
= \sigma_{j_{q}} \sigma_{j_{p}}$. Hence we have $\vee Y =
\sigma_{j_{1}} \sigma_{j_{2}} \ldots  \sigma_{j_{l}}$ in the weak order. 
It is easy to see that $\sigma_{j_{1}} \sigma_{j_{2}} \ldots  \sigma_{j_{l}}$ 
is a 321 avoiding permutation. 
Hence we have $\sigma_{j_{1}} \sigma_{j_{2}} \ldots  \sigma_{j_{l}} \in 
B^{\prime}_{n}$. \ \ \ \ \ \ \ $\Box$

\end{lem}

Next we will determine the NBB bases with respect to $\lhd$.

\begin{thm}\label{321-main-prop}

Put $X := \{ \sigma_{i_{1}}, \ldots , \sigma_{i_{k}} \} \subset A(B_{n})$ 
where $1 \le i_{1} < \cdots < i_{k} \le n-1$. We assume that 
$\vee X = \widehat{1}$.

Then we have  
$X$ is NBB $\Longleftrightarrow $
\begin{enumerate}
\item $i_{1} = 1$
\item $\{ i_{2}-1, i_{3} -1, \ldots i_{k}-1 \}$ is a sparse set of $[n-2]$.
\end{enumerate}

\end{thm}

{\bf Proof}

($\Longrightarrow$)

Because $X$ is an NBB base of $\widehat{1}$, we have $\sigma_{1} \in X$. 
Hence we have $i_{1} = 1$.
We put $Y := \{ \sigma_{i_{2}}, \sigma_{i_{3}} , \ldots 
\sigma_{i_{k}} \}$. 
If $\vee Y = \widehat{1}$ then we have $\sigma_{1} < \widehat{1}$ and 
$\sigma_{1} \not\in Y$. Then we have that $Y$ is not BB. 
This contradicts the assumption that $X$ is an NBB base. 
So we have $\vee Y \neq \widehat{1}$. 
By Lemma {\rmfamily \ref{321-avoiding-neighbor}} 
we have $i_{2} + 1 < i_{3}, \ i_{3} + 1 < i_{4}, \ \ldots 
i_{k-1} + 1 < i_{k} $. 
If $i_{2} \neq 2$ we have $\vee X \neq \widehat{1}$. 
Hence we have that the set 
$\{ i_{2}-1, i_{3} -1, \ldots i_{k}-1 \}$ 
is a sparse set of $[n-2]$.

($\Longleftarrow$)

Because $\sigma_{1}$ and $\sigma_{2}$ are elements of $X$, we have 
$\vee X = \widehat{1}$. 
Let $Y$ be a subset of $X$. 
We put $Y := \{  \sigma_{m_{1}}, \ldots \sigma_{m_{p}}  \} \subset X$ with 
$1 \le m_{1} < m_{2} < \ldots < m_{p} \le n-1$.
We have to show that $Y$ is not a BB base. 
If $m_{1} = 1$ it is clear. If $m_{1} \neq 1$ we have  
$m_{1} + 1 < m_{2}, \ m_{2} + 1 < m_{3}, \ldots m_{p-1} + 1 < m_{p}$. 
By Lemma {\rmfamily \ref{321-avoiding-sparse}} we have 
$\vee Y =  \sigma_{m_{1}} \sigma_{m_{2}} \ldots \sigma_{m_{p}}$. 
Hence we have $\{ x | \ x \in A(B_{n}), \ x \le \vee Y \} = 
\{ m_{1} , m_{2} , \ldots , m_{p} \}$. Hence we have that $Y$ is not a 
BB base. This completes the proof of our result. \ \ \ \ \ $\Box$

From Theorem {\rmfamily \ref{blass-sagan-main}} and Theorem 
{\rmfamily \ref{321-main-prop}} we obtain the following result.

\begin{thm}

We have 

\begin{eqnarray}
\mu(B_{n}) = \sum_{X \subset [n-2] \ X: \ {\rm sparse \ set}}
(-1)^{|X|+1} = (-1) F_{n-2}(-1).
\end{eqnarray}

\end{thm}

\begin{nota}

Note that 
For each $n \in \mathbb{N}$ the Tamari lattice $T_{n}$ is the poset of 
132 avoiding permutations with weak order on $S_{n}$. It is well known 
that for each $\sigma \in S_{3} \setminus \{ 123, 321\}$ the poset of 
$\sigma$ avoiding permutations is also the Tamari lattice.

\end{nota}

{\large \textbf{Acknowledgement}}

The author wishes to thank 
Professor Jun Morita for his 
valuable advices.

\renewcommand{\refname}{REFERENCE}

\end{document}